\documentclass{amsart}
\usepackage{amssymb}
\usepackage{amscd}
\swapnumbers
\newtheorem*{theorem*}{Theorem}
\newtheorem{theorem}{Theorem}[section]

\newtheorem{proposition}[theorem]{Proposition}

\theoremstyle{definition}
\newtheorem{definition}[theorem]{Definition}

\theoremstyle{remark}

\numberwithin{equation}{section}

\newcommand{\thref}[1]{Theorem~{\rm\ref{#1}}}

\newcommand{\deref}[1]{Definition~{\rm\ref{#1}}}

\newcommand{\seref}[1]{Section~{\rm\ref{#1}}}

\newcommand{\st}{\; | \;}                               

\newcommand{\ov}{\overline}
\newcommand{\del}{\partial}

\newcommand{\isoto}{\xrightarrow{\sim}}       

\newcommand{\CC}{\mathcal C}       
\newcommand{\E}{\mathcal E}       
\newcommand{\F}{\mathcal F}       
\newcommand{\C}{\mathbb{C}}       
\newcommand{\Z}{\mathbb{Z}}       
\newcommand{\R}{\mathbb{R}}       
\renewcommand{\O}{\mathcal{O}}   
\renewcommand{\P}{\mathcal{P}}

\newcommand{\ga}{\gamma}
\newcommand{\Ga}{\Gamma}
\newcommand{\de}{\delta}
\newcommand{\De}{\Delta}

\newcommand{\la}{\lambda}

\newcommand{\Si}{\Sigma}
\newcommand{\om}{\omega}

\newcommand{\eps}{\varepsilon}

\newcommand{\Mgn}{\mathcal M_{g,n}}

\newcommand{\MGgn}{\mathcal M_{g,n}^G}

\DeclareMathOperator{\Hom}{Hom}
\DeclareMathOperator{\Mor}{Mor}

\DeclareMathOperator{\Ad}{Ad}
\DeclareMathOperator{\End}{End}

\DeclareMathOperator{\Irr}{Irr}

\begin{document}

\title{On the modular functor associated with a finite group}

\author{Alexander Kirillov, Jr.}
   \address{Department of Mathematics, SUNY at Stony Brook, 
            Stony Brook, NY 11794, USA}
    \email{kirillov@math.sunysb.edu}
    \urladdr{http://www.math.sunysb.edu/\textasciitilde kirillov/}

\date{\today}

\maketitle

\section{Introduction}\label{s:introduction}
In this note, we discuss the complex-algebraic approach to the
modular functor constructed from a finite group $G$. This can be
considered as a ``baby'' case of a more interesting situation, modular
functors related to orbifold conformal field theories, which we intend
to pursue in a subsequent papers. 

We start by recalling some basic facts about the modular functor;
detailed exposition can be found, e.g., in \cite{BK}. Let $\CC$ be a
semisimple abelian category. Then the following structures are
essentially equivalent: 
\begin{enumerate}
\item Topological 2-dimensional modular functor, i.e. an assignment 
 \begin{equation}\label{e:mf-top1}
(\Si,p_i, V_i)\mapsto W(\Si, p_i, V_i).
\end{equation}
 Here $\Si$ is an oriented 2-dimensional surface with boundary, with
 marked points $p_i$ on each boundary circle $(\del\Si)_i$ and an
 object $V_i\in\CC$ assigned to $(\del\Si)_i$, and $W(\Si, p_i, V_i)$
 is a finite-dimensional complex vector space. This assignment should
 satisfy a number of properties, most important being functoriality
 and gluing axiom. 
 
\item Complex-algebraic modular functor, i.e. a collection of vector
  bundles with a projectively flat connection $W(C, V_i)$ on the
  moduli space $\Mgn$ of stable curves $C$ with marked points $p_i$
  and non-zero tangent vector $v_i\in T_{p_i}C$.  This assignment
  should satisfy a number of properties, most important being
  functoriality and factorization properties, which describes behavior of
  the connection near the boundary of Deligne--Mumford
  compactification of the moduli space.
\item A structure of a modular tensor category on $\CC$.
\end{enumerate}

Note that even though structures of a topological MF and a
complex-algebraic MF are equivalent, there is in general no simple way
to construct complex-analytic MF from a topological one. This is
essentially equivalent to constructing, for a given representation of
$\pi_1(\Mgn)$, a vector bundle with a flat connection with regular
singularities whose monodromy is described by this representation.
While the Riemann--Hilbert correspondence shows that such a local
system exists, it does not give a natural  construction of it.

A simplest example of a modular functor is a MF associated with a
finite group $G$ and a cohomology class $\om\in H^3(G, S^1)$. In
topological setting, it arises from the Chern--Simons theory with a
finite gauge group $G$; a detailed description can be found in
\cite{FQ} (we will review it below). In the language of modular
categories, the corresponding modular category is the category of
representations of the twisted Drinfeld double $D^{\om}(G)$ of the
finite group $G$ (see \cite{DPR}).  

In this note, we complete the picture by providing a complex-analytic
counterpart of the same modular functor. For simplicity, we only
describe untwisted case, i.e. $\om=1$. A key ingredient of the
construction is the moduli space of ``admissible $G$-covers'',
introduced in \cite{JKK}.

\section{Drinfeld double of a finite group}\label{s:double}
In this section, we briefly recall some facts about Drinfeld double of
the finite group. Throughout this section, $G$ is a finite group.

For an element $g\in G$, we denote by $\Ad_g\colon G\to G$ the adjoint
action of $g$ on $G$: $\Ad_g(x)=gxg^{-1}$. These operators naturally
form a groupoid, with the set of objects $G$ and $\Mor(x,y)=\{g\in
G\st Ad_g(x)=y\}$. We will denote this groupoid by $\Ad(G)$ (this
notation is not standard).

An action of the groupoid $\Ad(G)$ on a set $X$ is a decomposition 
$X=\sqcup_{g\in G}X_g$ and an action of $G$ on $X$ such that
$gX_h\subset X_{ghg^{-1}}$. Similarly, a {\em representation} of the
groupoid $\Ad(G)$ is a vector space $V$ with decomposition
$V=\bigoplus_{g\in G} V_g$ and a linear action of $G$ such that
$gV_h\subset V_{ghg^{-1}}$; in other words, a representation of
$\Ad(G)$ is the same as $G$-equivariant vector bundle on $G$.  

Similar to the usual construction for groups, we can define ``group
algebra'' of the groupoid $\Ad(G)$ as an algebra of formal linear
combinations of morphisms in $\Ad(G)$; we define product to be zero if
the morphisms are not composable. It is immediate that a
representation of $\Ad(G)$ is the same as a representation of the
group algebra (as an associative algebra with unit). It is also easy
to check that this group algebra is in fact isomorphic to the
semidirect product $D(G)=\C[G]\ltimes \F(G)$, where $\C[G]$ is the
group algebra of $G$ and $\F(G)$ is the algebra of functions on $G$.
Denoting by $\de_g\in \F(G)$ the delta function at $g$, the product in
$D(G)$ is given by $g\de_h=\de_{Ad_g(h)} g$. In fact, $D(G)$ has a
natural structure of Hopf algebra (see, e.g., \cite{BK}). This Hopf
algebra is called {\em Drinfeld double} of the group $G$.

As for any Hopf algebra, we can define the notion of invariants: if
$V$ is a representation of $D(G)$, then 
\begin{equation}\label{e:inv}
V^{D(G)}=\Hom_{D(G)}(\C, V)=(V_1)^{G}.
\end{equation}

It is well known that $D(G)$ is a semisimple associative algebra; thus
it has finitely many irreducible representations (up to isomorphism).
We will denote by $\Irr(D(G))$ the set of isomorphism classes of
irreducible representations. Explicit description of these
representations can be found, e.g., in \cite{BK}).  Semisimplicity
also implies that as a an algebra,
$$
D(G)\simeq \bigoplus_{\la\in \Irr(D(G))}\End(\rho_\la). 
$$

\section{Modular functor associated with a finite group: topological
  description}\label{s:topological} 
In this section, we briefly recall the definition of topological
modular functor associated with a finite group $G$, following
\cite{FQ} (with minor changes). 

Let $X$ be an oriented compact surface with boundary, and with a
marked point $p_i\in (\del X)_i$ in each boundary component of $X$. We
will call such a structure {\em marked surface} and will denote
$X^m=(X, \{p_i\})$. For such a surface, we denote by $\P(X^m)$ the
category with objects $P^m=(P, \{q_i\})$, where $P$ is a principal
$G$-bundle $\pi\colon P\to X$, and $q_i\in\pi^{-1}(p_i)$ is a chosen
lifting of the marked points $p_i$ to $P$. The morphisms in this
category are isomorphisms of $G$-bundles which map marked points to
marked points. Note that every such morphism is invertible, so $\P(X^m)$
is a groupoid, and that if each connected component of $X^m$ has at
least one boundary component, then objects of $\P(X^m)$ have no
non-trivial automorphisms, so $\P(X^m)$ is essentially a set. 

For a principal $G$-bundle with marked points $P^m=(P,\{q_i\})$, we
can define {\em monodromy} of $P^m$ around $(\del X)_i$ as follows.
Orientation of $X$ defines a natural direction of the boundary circle
$(\del X)_i$. Choose a parametrizattion $\ga\colon \R/2\pi\Z \to (\del
X)_i$ such that $\ga(0)=p_i$; then there is a unique lifting of $\ga$
to a map $\tilde\ga\colon [0,2\pi]\to P$ such that $\tilde\ga(0)=q_i$. We
define monodromy $m_i(P^m)$ by 
\begin{equation}
\label{e:monodromy}
\ga(2\pi)= m_i(P^m) \ga(0).
\end{equation}

Also, given $P^m$ and an element $g\in G$, we can define a new
marked $G$-bundle by using $g$ to change the marked point $q_i$:
\begin{equation}\label{e:rho}
\rho_i(g)P^m=(P, \{q_1,\dots, gq_i,\dots, q_n\}).
\end{equation}
Note that $\rho_i$ is not a morphism in the category $\P(X)$ (in
general, $\rho_i(g)P^m$ is not isomorphic to $P^m$) but a
functor $\P(X^m)\to \P(X^m)$. It is also immediate from direct
computation that
\begin{equation}\label{e:relation2}
m_i(\rho_i(g)P^m)= g\cdot m_i(P^m)\cdot g^{-1}.
\end{equation}

We denote by $\ov{\P(X^m)}$ the set of isomorphism classes in
$\P(X^m)$.  Then one easily sees that $m_i$ and $\rho_i$ descend to
$\ov{\P(X^m)}$, giving maps $\ov{\P(X^m)}\to G$ and action of $G$ on
$\ov{\P(X^m)}$. It follows from \eqref{e:relation2} that $m_i,\rho_i$
define an action of the groupoid $\Ad(G)$ on the set $\ov{\P(X^m)}$.

Now, let 
\begin{equation}\label{e:E(X)}
E(X^m)=\F(\ov{\P(X^m)}),
\end{equation}
where $\F(S)$ stands for the space of functions on $S$. 
 
Then for each boundary component $(\del X)_i$, the action of $\Ad(G)$
on $\ov{\P(X^m)}$ by $m_i,\rho_i$ defines on $E(X)$ a structure of
representation of $\Ad(G)$ and thus, by results of \seref{s:double},
of a representation of the algebra $D(G)$. We will denote this
representation by $\rho_i$. It can be written explicitly as follows:
\begin{equation}\label{e:rho2}
\begin{aligned}
(\rho_i(\de_h) f)P^m&=\de_{h, m_i(P^m)}f(P^m)\\
(\rho_i(g) f)(P^m)&=f(\rho_i(g^{-1})P^m).
\end{aligned}
\end{equation}

Let $\De(X)$ be the set of boundary components of $X$; for each
boundary component $(\del X)_i, i\in \De(X)$, consider a copy $D_i(G)$
of $D(G)$ and let 
\begin{equation}\label{e:D_delta}
D_{\De(X)}=\bigotimes_{i\in \De(X)}D_i(G).
\end{equation}
Taking tensor product of actions $\rho_i$ defined by \eqref{e:rho2},
we see that $E(X^m)$ has a natural structure of a
$D_{\De(X)}(G)$-module. 

\begin{definition}\label{d:topmf}
  For a marked surface $X^m$ and representations $V_i$ of $D(G)$
  assigned to boundary components $(\del X)_i$, we define the vector
  space
\begin{equation}\label{e:topmf}
W(X^m, \{V_i\})=\Hom_{D_{\De(X)}(G)}(E(X),\bigotimes V_i)
\end{equation} 
\end{definition} 

It immediately follows from the definition that 
\begin{equation}\label{e:decomposition}
E(X^m)=\bigoplus_{\la_1,\dots, \la_n}
 \rho^*_{\la_1}\otimes\dots\otimes \rho^*_{\la_n}
\otimes W(X^m, \rho_{\la_1},\dots, \rho_{\la_n})
\end{equation}
where each $\la_i$ runs over the set $\Irr(D(G))$ of isomorphims
classes of irreducible representations of $D(G)$, $\rho_{\la_i}$ is
the corresponding representation, and $\rho_{\la_i}^*$ is the dual
representation.

Finally, we define gluing. Let $X^m$ be a marked surface, and let
$c\subset X$ be a simple closed curve (``cut'') with a marked point
$p$. Cutting $X$ along $c$ gives a new surface $X^m_{cut}$, with $\del
X_{cut}=\del X\sqcup c'\sqcup c''$ and marked points $p', p''$ on $c',
c''$ respectively. Let $\P_{c',c''}(X^m_{cut})\subset \P(X^m_{cut})$
be the category of marked $G$-bundles such that
$m_{c'}(P)m_{c''}(P)=1$. We have an action of $G$ by functors on the
subcategory $\P_{c',c''}(X^m_{cut})$ given by
$\rho_{c',c''}(g)=\rho_{c'}(g)\rho_{c''}(g)$. Passing to the set of
isomorphism classes $\ov{\P(X^m_{cut})}$, we get a subset
$$
\ov{\P_{c',c''}(X^m_{cut})}\subset \ov{\P(X^m_{cut})}
$$
with an action of $G$. 

Any principal $G$-bundle on $X$ can be restricted to
$X_{cut}$. Analyzing which bundles on $X_{cut}$ can be obtained in
this way and taking into account marked points, one easily gets the
following proposition.

\begin{proposition}\label{p:gluing1}
Restriction gives a bijection 
$$
\ov{\P(X^m)}\isoto \ov{\P_{c',c''}(X^m_{cut})}/G
$$
\end{proposition} 

Passing to functions, we get the following result:
\begin{theorem}\label{t:gluing2}
We have a natural isomorphism
$$
E(X^m)\isoto (E(X^m_{cut}))^{D(G)}
$$
where the  action of $D(G)$ is given by $\rho_{c'}\otimes \rho_{c''}$. 
\end{theorem}

Note that $D(G)$ is not cocommutative, so the action of $D(G)$ on
$E(X^m_{cut})$ depends on the ordering of the cuts $c', c''$. However,
it is easy to see that the space of invariants does not depend on this
choice.

Finally, decomposing $E(X^m)$ into direct sum of irreducibles as in
\eqref{e:decomposition}, we immediately get the following result:
\begin{theorem}\label{t:gluing3}
One has a natural isomorphism of vector spaces 
$$
W(X^m, V_1,\dots, V_n)\simeq\bigoplus_\la W(X^m_{cut}, V_1,\dots, V_n,
\rho_\la, \rho_\la^*).
$$
\end{theorem}

Now it is easy to check the following result (see \cite{FQ} for
details):
\begin{theorem}
The assignment $X^m, V_1,\dots, V_n\mapsto W(X^m, \{V_i\})$, with
the gluing defined in \thref{t:gluing3}, satisfies all axioms of a
modular functor. The corresponding modular structure on the category of
representations of $D(G)$ coincides with the modular structure defined
by the quasi-triangular Hopf algebra structure on $D(G)$. 
\end{theorem}

\section{Moduli space of $G$-covers}\label{s:gcovers}
In this section, we briefly recall the definition of the moduli space
of stable $G$-covers, following \cite{JKK}, with one modification.
Namely, \cite{JKK} describes pointed curves (i.e., curves with marked
points $p_1,\dots, p_n$); we will also require choice of a non-zero
tangent vector $v_i\in T_{p_i}C$ at each marked point, or,
equivalently, 1-jet of local parameter at $p_i$. One can easily check
that all arguments  of \cite{JKK} apply in this case with obvious
changes, except for definition of gluing maps, which require more
extensive changes and which will be discussed in forthcoming papers.

First, let us recall the usual definition of marked curve and the
correpsonding moduli space. 

\begin{definition}
  A {\em marked curve} is a non-singular complex curve $C$ with
  distinct marked points $p_i$ and a choice of 1-jet of local
  parameter $dz_i$ at $p_i$. A marked curve is {\em stable} if the
  group of automorphisms preserving $p_i, dz_i$ is finite. 
\end{definition}

Note that we do not require that $C$ be connected. For connected
curves, stability means that the genus $g$ and number of marked points
$n$ are subject to restriction $(g,n)\ne (0,0), (0,1), (1,0)$.

We will denote by $\Mgn$ the moduli stack of connected genus $g$
marked curves with $n$ marked points. It is well-known that $\Mgn$ is
a smooth Deligne--Mumford stack.

Now let $G$ be a finite group. The following definition is a minor
modification of the one given in \cite{JKK}.

\begin{definition}\label{d:G-cover}
A non-singular marked $G$-cover is the following collection of data:
\begin{itemize}
\item A non-singular stable marked curve $C, \{p_i\}, \{dz_i\}$
\item A finite cover $\pi \colon \tilde C\to C$ and an action of $G$
  on $\tilde C$ which preserves projection $\pi$ and satisfies:
  \begin{itemize}
  \item Over $C-\{p_1,\dots, p_n\}$, $\pi$ is a principal $G$-bundle. 
  \item Near each $q_i\in \pi^{-1}(p_i)$, the map $\pi$ is locally
  analytically equivalent to $\C\to\C\colon \tilde z\mapsto z=\tilde
  z^{r_i}$ (the number $r_i$ is called the branching index at $p_i$).
 \end{itemize}
\item For each $i$, a choice of marked point $\tilde
  p_i\in\pi^{-1}(p_i)$ and 1-jet of local parameter $d\tilde z_i$ at
  $\tilde p_i$ such that $\tilde z^{r_i}=z_i$. 
\end{itemize}
\end{definition}

For future use, we note here that the tangent spaces $T_{p_i}C$ and
$T_{\tilde p_i}\tilde C$ are related by 
\begin{equation}\label{e:tangentspace}
\begin{aligned}
T_{p_i}C&=(T_{\tilde p_i}\tilde C)^{\otimes
  r_i}=(T_{\pi^{-1}(p_i)}\tilde C)/G\\
&T_{\pi^{-1}(p)}\tilde C =
\bigoplus_{q\in \pi^{-1}(p)} T_q\tilde C.
\end{aligned}
\end{equation}

As in the topological picture, given a $G$-cover $\tilde C\to C$ and a
marked point $p_i\in C$, one can define the monodromy $m_i(\tilde
C)\in G$ and action $\rho_i$  of $G$ by $\tilde p_i\mapsto g \tilde
p_i, \tilde z_i \mapsto \tilde z_i\circ g^{-1}$. Again, trivial check
shows that they satisfy relation $m_i(\rho_i(g)\tilde C)=g \cdot
m_i(\tilde C)\cdot g^{-1}$; thus, they define an action of the
groupoid $\Ad (G)$ on the category of $G$-covers.

The notion of $G$-cover can be easily defined for families of curves
(see \cite{JKK}). This allows one to define the moduli space of marked
$G$-covers. We will denote the moduli space of connected genus $g$
marked curves with $n$ marked points by  $\MGgn$. The same argument as in
\cite{JKK} shows that $\MGgn$ is a smooth Deligne--Mumford stack, and
the natural forgetting map defined by 
$$
st (\tilde C\to C)=C
$$ 
gives a morphism of stacks 
\begin{equation}\label{e:st}
st\colon\MGgn\to \Mgn
\end{equation}
which is a finite cover.

There is a direct relation of the complex-analytic picture
with the topological picture. Let $C, \{p_i\}, \{dz_i\}$ be a marked
curve.  Choose a local parameter $z_i$ at each $p_i$ with given 1-jet,
and a small enough positive real number $\eps$ such that $z_i$ is a
biholomorphic map $D_i\isoto \{z\in \C\st |z|<\eps\}$, for some
neighbourhoods $D_i$ of $p_i$, such that $D_i\cap
D_j=\varnothing$. Let $C^\circ=C\setminus \cup D_i$. This is an
oriented surface with boundary. Moreover, we have marked points on the
boundary specified by the condition $z_i=\eps$. Thus, $C^\circ$ is a
marked surface. 

Note that $C^\circ$ depends on the choice of local parameters $z_i$
and $\eps$, but it can be shown that different choices produce marked
surfaces which are isomorphic, and isomorphism is canonical up to
homotopy (see \cite{BK}). In particular, this implies that the modular
functor space $W(C^\circ, V_1,\dots, V_n)$ is canonically defined.

\begin{theorem}
  Let $C, C^\circ$ be as above. Then the category of admissible marked
  $G$-covers $\tilde C\to C$ is naturally equivalent to the category
  $\P(C^\circ)$ of marked $G$-bundles over $C^\circ$. 
\end{theorem}
\begin{proof}
  First, restrictings a $G$-cover $\tilde C\to C$ to $C^\circ\subset C$
  and forgetting the complex structure  defines a functor from
  $G$-covers to $\P(C^\circ)$. One easily sees that conversely, given
  a principal $G$-bundle on $C^\circ$, there is a unique way to extend
  it to a (topological) branched cover $\tilde C$, and then there is a
  unique complex structure on $\tilde C$ compatible with projection
  $\tilde C\to C$. 
\end{proof}

\section{Modular functor associated with a finite group:
  complex-algebraic description}\label{s:complex} 
In this section, we finally formulate the main result of this note,
namely a definition of complex-analytic modular functor associated
with a finite group $G$.

Let $\MGgn$ be the moduli space of $G$-covers defined in
\seref{s:gcovers}. Consider the structure sheaf $\O$ on $\MGgn$; it is
obviously a module of the sheaf $D_{\MGgn}$ of differentail operators
on $\MGgn$. Define a sheaf $\E$ of $D$-modules on $\Mgn$ by
\begin{equation}\label{e:E}
\E=st_*(\O)
\end{equation}
where $st\colon \MGgn\to \Mgn$ is the forgetting map
\eqref{e:st}.

Since $st\colon\MGgn\to \Mgn$ is a finite cover, one easily sees that
$\E$ is a lisse $D$-module, i.e. a sheaf of sections of some vector
bundle $E$ with flat connection. Action of the groupoid $\Ad (G)$ on the
moduli space $\MGgn$ gives an action of $\Ad (G)$ (and thus, of the
algebra $D(G)$) on the sheaf $\E$. Define $\De(C)$ be the set of
marked points of $C$ and define, in analogy with \eqref{e:topmf}, the
modular functor sheaf
\begin{equation}\label{e:E2}
\E(V_1, \dots, V_n)=\Hom_{D_{\De(C)}(G)}(\E,\O\otimes
V_1\otimes\dots\otimes V_n). 
\end{equation}

This sheaf is naturally the sheaf of holomorphic sections of the vector
bundle $E$ with fiber 
\begin{equation}\label{e:E2}
E_C(V_1, \dots, V_n)=\Hom_{D_{\De(C)}(G)}(E, V_1\otimes\dots\otimes V_n). 
\end{equation}

\begin{theorem}\label{t:comp=top}
\begin{enumerate}
\item Let $C$ be a marked curve. Then one has canonical isomorphisms
  of $D_{\De(C)}(G)$-modules 
\begin{equation}\label{e:comp=top1}
E_C=E(C^\circ)
\end{equation}
\item Let $C$ be a marked curve. Then one has canonical isomorphisms
  of vector spaces
\begin{equation}\label{e:comp=top2}
E_C(V_1,\dots, V_n)=W(C^\circ,V_1,\dots, V_n )
\end{equation}
\item Under isomorphism \eqref{e:comp=top2}, the representation of 
  the mapping class groupoid $\Ga$ given by
  monodromy of the local system $E_C(V_1,\dots, V_n)$  is identified
  with the representation given by the modular functor
  $W(C^\circ,V_1,\dots, V_n )$. 
\end{enumerate}
\end{theorem}

In short, this theorem states that under the correspondence between
topological MF and complex-analytic MF described in \cite{BK}, the
modular functor defined by \eqref{e:E} corresponds to the finite group
modular functor defined in \deref{d:topmf}. 

Note that this theorem does not address the question of defining
the gluing isomorphism in the complex-analytic approach. This will be
discussed in subsequent papers.

\end{document}